\documentclass[myspacing,11pt]{ectaart1}

\RequirePackage[OT1]{fontenc} \RequirePackage{amsthm}
\RequirePackage[cmex10]{amsmath} \RequirePackage{natbib}
\RequirePackage[colorlinks]{hyperref} \RequirePackage{hypernat}
\RequirePackage{amsfonts} \RequirePackage{amssymb} \RequirePackage{a4wide}
\RequirePackage{graphicx}

\startlocaldefs \numberwithin{equation}{section}
\theoremstyle{plain}
\newtheorem{definition}{Definition}[section]
\newtheorem{theorem}{Theorem}[section]

\newtheorem{remark}{Remark}[section]
\newtheorem{lemma}{Lemma}[section]

\endlocaldefs

\def\@bysame#1{\vrule height 1.5pt depth -1pt width 3em \hskip
0.5em\relax}

\newcommand{\N}{ \mathbb{N} }

\newcommand{\Z}{ \mathbb{Z} }

\newcommand{\R}{ \mathbb{R} }

\newcommand{\trunc}[1]{ {\lfloor #1 \rfloor} }
\newcommand{\wh}[1]{ \widehat{ #1 } }
\newcommand{\wt}[1]{ \widetilde{ #1 } }

\newcommand{\calF}{\mathcal{F}}

\newcommand{\calH}{\mathcal{H}}

\newcommand{\calI}{\mathcal{I}}

\newcommand{\calS}{\mathcal{S}}

\newcommand{\eins}{{ 1}}

\newcommand{\Cov}{{\mbox{Cov\,}}}

\newcommand{\argmin}{ \operatorname{argmin} }

\endlocaldefs

\begin{document}

\begin{frontmatter}

\title{Sequential Data-Adaptive Bandwidth Selection by Cross-Validation for Nonparametric Prediction}
\runtitle{Sequential Bandwidth Selection}


\begin{aug}
\author{\fnms{Ansgar} \snm{Steland}
\ead[label=e1]{steland@stochastik.rwth-aachen.de}}
\ead[label=u1,url]{www.stochastik.rwth-aachen.de}

\runauthor{A. Steland}

\address{Institute of Statistics\\ RWTH Aachen University \\ W\"ullnerstr. 3, D-52056 Aachen, Germany\\
\printead{e1}} \printaddresses \bigskip

\end{aug}


\runauthor{A. Steland}

\begin{abstract}
We consider the problem of bandwidth selection by cross-validation from a sequential point of view in a nonparametric regression
model. Having in mind that in applications one often aims at estimation, prediction and change detection simultaneously,
we investigate that approach for sequential kernel smoothers in order to base these tasks on a single statistic.
We provide uniform weak laws of large numbers and weak consistency results for the cross-validated bandwidth.
Extensions to weakly dependent error terms are discussed as well. The errors may be $ \alpha $-mixing or $L_2$-near
epoch dependent, which guarantees that the uniform convergence of the cross validation sum and the consistency of the
cross-validated bandwidth hold true for a large class of time series. The method is illustrated by analyzing photovoltaic data.
\end{abstract}




\end{frontmatter}

\section{Introduction}

The nonparametric regression model, often estimated by estimators of the Nadaraya-Watson type, forms an attractive framework for diverse areas such as engineering, econometrics, environmetrics, social sciences and biometrics. The present paper is devoted to a detailed study of a sequential bandwidth selector for kernel-weighted sequential smoothers related to the Nadaraya-Watson estimator. 
However, there are some subtle differences compared to the treatment of that estimator in nonparametric regression, since our Nadaraya-Watson type statistic is a prediction statistic, which we use to detect a change in the mean of the observations. Addressing the detection problem, we consider a setup which differs from that used  in classic nonparametric regression; especially, our setup leads to bandwidth choices not approaching $0$, as the sample size increases. Thus, although there is an interesting and obvious link to the classic regression problem, which we shall discuss in the next paragraph, the asymptotic results as well as the bandwidth selection problem are different and new.

Let us assume that observations $ Y_n = Y_{Tn} $, $ 1 \le n \le T $, arrive sequentially until the maximum sample size $ T $ is reached and satisfy the model equation
\begin{equation}
\label{ModelPaper}
  Y_n = m(x_n) + \epsilon_n, \qquad n = 1, 2, \ldots, T, \ T \ge 1,
\end{equation}
for the fixed design
\[
  x_n = x_{Tn} = G^{-1}( n/T ), \qquad 1 \le n \le T,
\]
induced by some design distribution function $ G $ and some function $ m:[0,\infty) \to \mathbb{R} $ (assumptions on $m$ will be given below). Tentatively, we make the assumption that the errors  $ \{ \epsilon_n : n \in \mathbb{N} \} $ form a sequence of i.i.d.($F$) random variables such that $ E( \epsilon_n ) = 0 $. We shall provide general results for weakly dependent time series, namely for strong mixing as well as near epoch dependent (NED) processes, but intend to postpone this issue to the end of Section~\ref{Sec:Asymptotics}, in order to focus on the idea of sequential cross-validation first.
 
Notice that in many applications the design points are either given or selected according to some external optimality criterion such that it is not restrictive to assume that $G$ is known. For instance, in econometrics the time instants where prices are quoted are usually fixed and known. Similarly, when discretizing signals or logging 
internet traffic, the variables of interest are sampled at known time points. In other applications, it may be preferable to use more design points in regions where $m$ is expected to be more volatile than in other regions, or in regions where higher accuracy is required. The latter issue may matter, for instance, when analyzing the nonlinear relationship between a medical response variable and an explanatory variable such as age or blood pressure, or, in social sciences, e.g. when studying the influence of the duration of unemployment on variables measuring quantities such as political opinion or social networking. Hence, we can and shall assume that $ x_n = n/T $, otherwise replace $m$ by $ \widetilde{m} = m \circ G^{-1} $, and interpret the regressor as a time variable. 

Clearly,  $ m(t) $ models the process mean (signal) of the underlying observations.
In practice, an analysis has often to solve three problems. (i) Estimation of the current process mean. (ii) One-step prediction of the process mean. (iii) Signaling when there is evidence that the process mean differs from an assumed (null) model. Usually, different statistics are used for those problems. For nonparametric estimation various methods have been studied including kernel estimators, local polynomials, smoothing splines and wavelets; we refer to
\citet{DonohoJohnston1994}, \citet{Eubank1988}, \citet{HaerdleAppliedNonpR1991} and \citet{WandJones1995}, amongst others.
Concerning procedures proposed in the literature to detect changes, there are various kinds of proposals. Some rely on closely related versions of those estimators, e.g. \citet{WuChu93}, \citet{MuellerStadtmueller99}, \citet{Steland2005JSPI} or \citet{Steland2009JNonpStat}, whereas other proposals construct special methods  as in \citet{PawlakRafaStelandJNonpStat04,PawRafSte08}, \citet{RafSte09} and \citet{PawRafSte10}, or apply classic control chart statistics of the CUSUM, MOSUM or EWMA type. For the latter approach we refer to \citet{HorvathKokoszka2002Statistics} and \citet{BrodskyDarkhovsky2000Book}, amongst many others. Frequently, change-point asymptotics can be based on the classic invariance principle of Donsker and its various generalizations to dependent time series, and functional asymptotics, which plays an important role in functional data analysis as well; we refer to \citet{Steland10a}, \citet{RafSteland10}, \citet{Bosq1998} and \citet{HorvHuskKok10}. There is also a rich literature on the estimation of regression functions that are smooth except some discontinuity (change-) points. See, for example, the recent work of \citet{GijbelsGoderniaux2004} or \citet{AntochGregoireHuskova2007}.

Separating change detection from estimation and prediction has benefits and drawbacks. Of course, it allows us to apply a detector which has certain optimality properties, but this requires knowledge of the model after the change, which is often too restrictive for practical applications. Further, sequentially analyzing two or even more sequences of statistics may be prohibitive in real world applications. Thus, to ease interpretation and applicability, the present paper investigates the idea to base a detector on a prediction statistic which can be used as an estimator as well. Our reasoning is that a method which fits the data well and has convincing prediction properties should also possess reasonable detection properties for a large class of alternatives models. 

The proposed kernel smoother requires to select a bandwidth parameter which controls the degree of smoothing. As well known, the bandwidth choice is crucial for performance. The topic has been extensively studied for the classic problem of nonparametric regression where the data gets dense as the sample size increases. Cross-validation belongs to the solutions which have been widely adopted by practitioners. To the best of our knowledge, sequential cross-validation as treated in the present paper has not yet been studied in the literature. We propose to select the bandwidth sequentially by minimizing a sequential version of the cross-validation criterion. In this way, one can update the bandwidth when new data arrive. Since we have in mind the detection of changes where consistent estimation is not really required, we base our analysis on a framework which is quite common in time series analysis and engineering signal analysis, but differs from the nonparametric regression setting used to obtain consistency results. This is motivated by the fact that in many applications the data are observed at a scale which does not converge to $0$, as the number of available observations approaches infinity.

The present paper aims at presenting first results on sequential cross-validation focusing on uniform consistency. We establish weak and $ L_2 $- consistency of the proposed sequential cross-validation criterion, uniformly over the time points where cross-validation is extended. Our results allow us to choose the number $N$ of time points as a function of the maximum sample size $T$ as well as to select their locations depending on $T$, as long as $N$  grows not too fast compared to $T$. We also extend the results to obtain weak consistency uniformly over compact sets for the bandwidth (parameter). The results yield a consistency result for the optimal bandwidth under quite general conditions on the above model. Finally, we extend consistency to $ \alpha $-mixing time series and near epoch dependent series.

The plan of the paper is as follows. Section~\ref{Sec:MethAss} discusses our assumptions and introduces in detail the sequential kernel smoother of interest. In Section~\ref{Sec:SeqCrossVal}, we introduce the sequential cross-validation approach. Our asymptotic results for i.i.d. errors are provided and discussed in Section~\ref{Sec:Asymptotics}. Section~\ref{Sec:Dependent} elaborates on extensions to dependent data. Those extensions work under very general assumptions, thus ensuring that the proposed method is valid for many real data sets. Detailed proofs of the main results are postponed to an appendix. Section~\ref{Sec:Illustration} discusses an application of the proposal to a case study from photovoltaic engineering dealing with power output measurements of photovoltaics modules (solar cells).

\section{Assumptions and Sequential Smoothers}
\label{Sec:MethAss}

Our mathematical framework is as follows.  Since the information about the problem of interest is often not sufficient to setup a (semi-) parametric model for the process mean $ m $ and the distribution of the error terms, which would allow us to use methods based on, e.g., likelihood ratios, a nonparametric framework is employed. We assume that model (\ref{ModelPaper}) holds true for a function $m$ with
\begin{equation}
\label{ST:AssumptionsM}
  m \in \text{Lip($[0,\infty);\R$)}, \quad \text{either $m>0$ or $m<0$}, \qquad \text{and $ \| m \|_\infty < \infty $ },
\end{equation}
where $ \text{Lip}(A,B) $, $ A, B \subset \R $, denotes the class of Lipschitz continuous functions $ A \to B $. Clearly, cross-validation is meaningless if $ m = 0 $. Having in mind possible applications where one aims at detecting quickly that the process level $m$ gets either too large or too small, we assume that either $ m > 0 $ or $ m < 0 $ and, w.l.o.g., confine ourselves to the case  $ m > 0 $ in what follows. Recalling that extensions to weakly dependent processes will be given in Section~\ref{Sec:Dependent}, let us assume at this point that $ \{ \epsilon_n \} $ are mean zero i.i.d. with common distribution function $ F $ satisfying
\begin{equation}
\label{ST:Moments}
  \int x^4 dF(x) < \infty.
\end{equation}
Under the general condition (\ref{ST:AssumptionsM}), one should use statistical methods which avoid (semi-) parametric specifications of the shape of $m$. Instead, nonparametric smoothers $ \wh{m}_n $ which estimate some monotone functional of the process mean and which are sensitive with respect to changes of the mean are of interest. 

Thus, given a kernel function $ K: [0, \infty) \to [0, \infty) $ and a bandwidth $ h > 0 $ the following sequential kernel smoother
\[
  \wt{m}_{i} = \wt{m}_{i,h} = \frac{1}{h} \sum_{j=1}^n K([j-i]/h) Y_i, \qquad i = 1, 2, \ldots
\]
and the associated normed version
\[
  \wh{m}_i = \wh{m}_{i,h} = \wt{m}_{h} \ \biggl/ \ \frac{1}{h} \sum_{j=1}^i K([j-i]/h),
\]
respectively, which are closely related to the classic Nadaraya-Watson estimator, are the starting points of our discussion.

\begin{remark}
At this point, it is worth noting that various classic control chart statistics are obtained as special cases. Denoting the target value by $ \mu_0 $, the CUSUM chart is based on $ C_i = \sum_{j=1}^i [X_j-(\mu_0+L)] $ where $ \{ X_n \} $ denotes the observed process and $ L $ is the reference value. This chart corresponds to the choice $K(z) = 1 $, $ z \in \R $, and $ Y_j = X_j - (\mu_0 + L) $ for all $j$. The EWMA recursion, $ \widehat{m}_i = \lambda Y_i + (1-\lambda) \widehat{m}_{i-1} $ with starting value $ \widehat{m}_0 = Y_0 $, $ \lambda \in (0,1) $ a smoothing parameter, corresponds to the kernel $ K(z) = e^{-|z|} $ and the bandwidth $ h = 1/\log(1-\lambda) $. By defining the weights by means of a kernel function, we get a rich class of statistics covering classic control statistics as special cases.
\end{remark}

The canonical one-sided detectors (stopping times) studied in \citet{SchmidSteland2000}, \citet{Steland2004} and \citet{Steland2005JSPI} have the form
\[
  S_T^- = \inf \{ \trunc{s_0 T} \le i \le T : \wh{m}_i > c \} \qquad \text{and} \qquad 
  S_T^+ = \inf \{ \trunc{s_0 T} \le i \le T : \wh{m}_i < c \},
\]
respectively. Here $c$ is a threshold (control limit), $ s_0 \in (0,1) $ determines through $ \trunc{Ts_0} $ the start of monitoring,  and $ \trunc{x} $ denotes the integer part (floor function) of $x$. Notice that $ S_T^{-} $ and $ S_T^+ $ are indeed stopping times, i.e., for instance, $ \{ S_T^- < n \} \in  \sigma(Y_1, \dots, Y_n) $ for all $ n \in \N $. A related stopping time is used in our illustration, cf. Section~\ref{Sec:Illustration}.

Our assumptions on the smoothing kernel are as  follows. 
\begin{equation}
\label{ST:AssumptionsK}
 \text{$K \in $ Lip($[0,\infty); [0,\infty)$), $ \| K \|_\infty < \infty $, $ \text{supp}(K) \subset [0,1] $, and $K>0$ on $(0,1)$.} 
\end{equation}
These assumptions are quite standard and satisfied by many kernels used in practice. Our result on the uniform weak law of large numbers for dependent time series even works under a weaker condition discussed there. It is well known that the choice of the bandwidth is of more concern than the choice of the kernel. However, in \citet{Steland2005JSPI} the problem of optimal kernel choice for detectors based on kernel-weighted averages has been studied in greater detail. In this work it is shown that the optimal kernel which minimizes the asymptotic normed delay depends on the alternative, i.e. on the mean of the process after the change. Particularly, CUSUM type procedures are not optimal in general models. Although the detection statistic studied there slightly differs from the prediction statistic studied in the present paper, those results may be used to some extent in order to select a kernel, if there is some a priori knowledge on possible models for the mean after the change. However, in what follows we assume that a kernel satisfying Assumption (\ref{ST:AssumptionsK}) has been selected, such that it remains to select a bandwidth.

For the bandwidth $ h > 0 $ we assume that
\begin{equation}
\label{ST:Bandwidth}
  |T/h - \xi | = O(1/T) 
\end{equation}
for some constant $ \xi \in [1, \infty ) $, where the $ O(1/T) $-requirement (instead of $ o(1) $) rules out artificial choices such as $ h = T/(\xi + T^{-\gamma}) $, $ \gamma > 0 $, leading to arbitrary slow convergence. It is worth discussing that assumption, which is rather different than the $h \to 0 $ such that $ nh \to \infty $ assumption encountered in nonparametric regression. In our setup, we work with an equidistant design where the distance between the time points does not converge to $0$, i.e., we use no {\em in-fill} asymptotics. Assumption (\ref{ST:Bandwidth}) now guarantees that the number of observations on which $ \wh{m}_T $ depends converges to $ \infty $, as $ T \to \infty $. In practice, one can select $ \xi $ and put $ h = T/\xi $. Notice that in our asymptotic setup the parameter $ \xi $ determines the degree of {\em localization} of the procedure. If one uses an (approximation to the) uniform kernel, $ \xi $ fixes the percentage of observations used in each step of the detection procedure.

The asymptotic distribution theory for procedures based on the sequential kernel smoother $ \widehat{m}_n $ has been
studied in \citet{Steland2004}, \citet{Steland2005JSPI} and \citet{Steland05JTSA}. Those results allow us to construct classic statistical hypothesis tests as well as monitoring procedures to detect changes in the process mean, such that certain statistical properties are asymptotically satisfied. Specifically, it is shown that for a large class of weakly dependent error processes $ \{ \epsilon_t \} $ the process $ \{ \sqrt{T}  \widehat{m}_{\trunc{Ts},h} : s \in [0,1] \} $ satisfies a functional central limit theorem with Gaussian limit when the underlying observations have mean $0$, i.e.,
\[
  \sqrt{T} \widehat{m}_{\trunc{Ts},h} \Rightarrow \mathbb{M}(s),
\]
as $ T \to \infty $, for some centered Gaussian process $ \{ \mathbb{M}(s) : s \in [0,1] \} $ which depends on $ \xi $; that result covers the limiting distribution of the classic Nadaraya-Watson estimator in our setting as a special case. However, it turns out that the asymptotic law and therefore the control limit ensuring that the asymptotic type I error rate satisfies $ \lim_{T \to \infty} P( S_T^- \le T ) = \alpha $ and $  \lim_{T \to \infty} P( S_T^+ \le T ) = \alpha $, respectively, depends on $ \xi $, where $ \alpha \in (0,1) $. The question arises, how one can or should select the bandwidth $ h \sim T $ and the parameter $ \xi $, respectively.

\section{Functional Sequential Cross-Validation}
\label{Sec:SeqCrossVal}

In the present paper, we propose to select the bandwidth  $ h > 0 $ such that the $ Y_t  $ are well approximated by  sequential predictions calculated from past data $ Y_1, \dots, Y_{t-1} $. For that purpose, we propose a  sequential version of the cross-validation criterion based on sequential leave-one-out estimates. 

The idea of cross-validation is to choose parameters such that the corresponding estimates provide a good fit on average. To  achieve this goal, one may consider the average squared distance between observations, $Y_i$, and predictions as  an approximation of the integrated squared distance. To avoid over-fitting and interpolation, the prediction of $ Y_i $ is determined using the reduced sample where  $ Y_i $ is omitted. Aiming at selecting the bandwidth $h$ to obtain a good fit when using  sequential prediction estimates, we consider 
\begin{equation}
\label{DefWHat}
  \widehat{m}_{h,-i} = N_{T,-i}^{-1} \frac{1}{ h } \sum_{j=1}^{i-1} K( [j-i]/h ) Y_j, \qquad i = 2, 3, \ldots
\end{equation}
with the constant $ N_{T,-i} = h^{-1} \sum_{j=1}^{i-1} K([j-i]/h) $. $ \widetilde{m}_{h,-i} $ is defined accordingly. Notice that these statistics are $ \sigma(Y_1, \dots, Y_{i-1}) $-measurable, i.e. adapted.

The statistic $ \widehat{m}_{h,-i} $ can be regarded as a sequential leave-one-out estimate. In (\ref{DefWHat}) we define the kernel weights using the bandwidth $h$; the kernel $K$ puts a weight on the distance between $ j/h $ and $i/h$. Due to assumption (\ref{ST:Bandwidth}), this is asymptotically the same as putting a weight on the distance between the time point $ i/T $ at which we want to predict the response and the time point $j/T$.
The corresponding detectors are given by
\[
  S_T^+ = \inf \{ \trunc{s_0 T} \le i \le T : \wh{m}_i > c \} \qquad \text{and} \qquad 
  S_T^-  = \inf \{ \trunc{s_0 T} \le i \le T : \wh{m}_i < c \},
\]
respectively. Given the predictions $ \widehat{m}_{h,-i}  $ we may define the {\em sequential leave-one-out cross-validation criterion} 
\[
  CV_s(h) = CV_{T,s}(h) = \frac{1}{T} \sum_{i=2}^{\trunc{Ts}} ( Y_i - \widehat{m}_{h,-i} )^2, \qquad s \in [s_0,1],  \quad h > 0.
\]
The cross-validation bandwidth at time $s$ is now obtained by minimizing $ CV_s(h) $ for fixed $ s $. To be precise, we are interested in the
following optimization problem where one minimizes over a set of arrays.
Let $ \calH_{s_0,\xi} $ be the family of all arrays $ \{ h_{Tn} : \trunc{s_0T} \le n \le T, \ T \ge 1 \} $ with
\[
  \lim_{T \to \infty} \frac{ T }{ h_{Tn} } = \xi \qquad \text{for some $ \xi > 0 $}.
\]
Now one considers minimizers $ \{ h_{Tn}^* \} \in \calH_{s_0,\xi} $ of the cross-validation criterion such that 
\[
  CV_{n/T}(h_{Tn}^*) \le CV_{n/T}( h_{Tn} ), \qquad \trunc{s_0T} \le n \le T, \ T \ge 1, 
\]
for all $ \{ h_{Tn} \} \in \calH_{s_0,\xi} $. That procedure yields the cross-validated bandwidth $ h^*_{T,\trunc{Ts}/T} $ for fixed $s$. Therefore,
\[
  h_T^*(s) = h^*_{T,\trunc{Ts}/T}, \qquad s \in [s_0,1],
\]
is our functional sequential estimate for the bandwidth.

The idea to proceed is now as follows.
We shall show that $ CV_{T,s}(h) $ converges to some function $ CV_\xi(s) $ which depends on $ \xi = \lim T/h $ provided that limit exists. Now we expect that under certain conditions  $ T/h_{T,\trunc{Ts}/{T}}^* $ converges to a minimizer of the function $ \xi \mapsto CV_\xi(s) $.  That minimizer yields the asymptotically optimal constant of proportionality $ \xi_s^* $.

\begin{remark} The following remarks are in order.
\begin{itemize}
\item[(i)]
Notice that $ CV_s(h) $ is a sequential unweighted version of the criterion studied by \citet{HaerdleMarron1985} in the classic regression function estimation framework. We do not consider a weighted CV sum, since we have in mind that the selected bandwidth is used to obtain a good fit for past and current observations. However, similar results as those presented here can be obtained for a weighted criterion such as $ T^{-1} \sum_{i=1}^{\trunc{Ts}} K( [\trunc{Ts}-i] /h) ( Y_i - \widehat{m}_{h,-i} )^2 $ as well. 
\item[(ii)]
At first glance, our approach is similar to one-sided cross-validation proposed by \citet{HartYi1998JASA} for bandwidth selection of nonparametric regression estimators in the classic regression framework. However, we are interested in sequential bandwidth selection and aim at studying the random function $ s \mapsto \argmin_h CV_s(h) $.
\end{itemize}
\end{remark}

Let us close this section with a discussion how to implement the approach in practice.
Cross-validation is expensive in terms of computational costs and minimizing $ C_{T,s} $ for all $s \in \{ n/T : \trunc{s_0T} \le n \le T \} $ is not feasible in many cases. Therefore and to simplify exposition, let us fix a finite number of time points $ s_1, \dots, s_N $ such that
\[ 
  0 < s_0 < s_1 < \cdots < s_N \le 1,
\] 
$ N \in \mathbb{N} $. However, for small $T$ a small value for $N$ is appropriate, whereas one would prefer a larger value for $N$ when $ T $ is large. Thus, it would be nice if $N$ could dependent on $T$. Indeed, we shall later relax this assumption and allow that $N$ is an increasing function of $ T $. At time $ s_i $ the cross-validation  criterion is minimized to select the bandwidth, $ h_i^* = h_i^*( Y_1, \dots, Y_{s_i} ) $, and that bandwidth is used during the  time interval $ [s_i,s_{i+1}) $, $ i = 1, \dots, N $.

\section{Asymptotic Results for I.I.D. Errors}
\label{Sec:Asymptotics}

The present section is devoted to a careful discussion of the asymptotic results of the present paper for i.i.d. errors. We provide several theorems on weak uniform consistency of the sequential cross-validation approach including results which allow that the number of time points where the cross-validated bandwidth is computed gets larger as $T$ increases. Further, we show that the optimal bandwidth behaves nicely in the limit in the sense described in the previous section by establishing an argmin consistency result which identifies the asymptotic constant of proportionality under certain regularity conditions. 

\subsection{Uniform convergence}

Notice that due to
\[
  CV_s(h) = \frac{1}{T} \sum_{i=1}^{\trunc{Ts}} Y_i^2 - \frac{2}{T} \sum_{i=2}^{\trunc{Ts}} Y_i \widehat{m}_{h,-i}
    + \frac{1}{T} \sum_{i=2}^{\trunc{Ts}} \widehat{m}_{h,-i}^2
\]
minimizing $ CV_s(h) $ is equivalent to minimizing
\[
  C_{T,s}(h) = - \frac{2}{T} \sum_{i=2}^{\trunc{Ts}} Y_i \widehat{m}_{h,-i}
    + \frac{1}{T} \sum_{i=2}^{\trunc{Ts}} \widehat{m}_{h,-i}^2.
\]
Thus, we will study $ C_{T,s}(h) $ in the sequel.

Our first result identifies the limit in mean of $ C_{T,s}(h) $ from which the asymptotically optimal constant of proportionality can be eventually calculated.

\begin{theorem}
\label{ST:Th0}
Assume (\ref{ST:AssumptionsM}) and (\ref{ST:AssumptionsK}). Then
\begin{align}
\label{St:DefC}
  E( C_{T,s}(h) ) &\to C_\xi(s) =
  - 2 \frac{ \int_0^s \int_0^r \xi K( \xi(r-u) ) m( \xi u ) \, du dr }
           { \int_0^s \xi K(\xi(s-r)) \, dr } \\
  \nonumber
  & + \frac{ \int_0^s \xi^2 \int_0^r \int_0^r K( \xi(r-u) ) K( \xi(r-v) ) m(u) m(v) \, du \, dv \, dr }
           { \int_0^s \xi K( \xi(s-r) ) \, dr },
\end{align}
as $ T \to \infty $, uniformly in $ s \in [s_0,1] $.
\end{theorem}

It is worth mentioning that point-wise convergence holds true under weaker conditions, e.g., if $ K$ is bounded and continuous and $ m $ is continuous with $ \int_0^1 m^2(t) \, dt < \infty $. Further, Theorem~\ref{ST:Th0} does not require independence as long as $ \{ \epsilon_n \} $ are pair-wise uncorrelated.
For an example illustrating the function $ C_\xi(s) $ we refer to our preliminary study \citet{Steland2009JNonpStat}. 

We will now study the (uniform) mean squared convergence of the random function $ C_{T,s}(h) $. 
Define $ \calS_N = \{ s_i : 1 \le i \le N \} $. $L_2$-consistency holds true at the usual rate.

\begin{theorem}
\label{ST:Th1}
Assume (\ref{ST:AssumptionsM}), (\ref{ST:Moments}) and (\ref{ST:AssumptionsK}). Then, for any fixed integer $ N $, we have the law of large numbers in $ L_2 $,
\[
  E \max_{s \in \calS_N} |C_{T,s}(h) - E(C_{T,s}(h))|^2 = O( T^{-1} ),
\]
as $ T \to \infty $.
\end{theorem}

The question arises whether we may increase the number of time points where cross-validation is conducted, if the maximum sample size $T$ increases. The following theorem provides such a uniform law of large numbers, but we no longer have a convergence rate.

\begin{theorem} 
\label{ST:Th2}
Assume $ N = N_T $ is an increasing function of $ T $ and
\begin{equation}
\label{ST:SN}
  0 < s_0 < s_{N1} < \cdots < s_{NN} \le 1, \qquad N \ge 1,
\end{equation}
and put $ \calS_N = \{ s_{Ni} : 1 \le i \le N \} $. Given Assumptions (\ref{ST:AssumptionsM}), (\ref{ST:Moments}) and (\ref{ST:AssumptionsK}), we have the uniform law of large numbers in $ L_2 $,
\[
  E \sup_{s \in \calS_N} |C_{T,s}(h) - E(C_{T,s}(h))|^2 = o(1),
\]
as $ T \to \infty $, provided 
$$
  \frac{N_T}{T} = o(1).
$$
\end{theorem}

\begin{remark}
It is worth mentioning that the location of the $N_T$ time points may depend on $N$, as long as they remain deterministic. If they are selected 
at random, the results remain valid a.s., as long as $ \calS_N $ and $ \{ \epsilon_t \} $ are independent, since then one can condition on $ \calS_N $.
\end{remark}

Combining the above statements, we obtain the following result.

\begin{theorem} 
\label{ST:Th3}
Suppose Assumptions (\ref{ST:AssumptionsM}), (\ref{ST:Moments}) and (\ref{ST:AssumptionsK}) hold true and, additionally, \eqref{ST:SN} is satisfied. Then
\[
  E \sup_{s \in \calS_N} |C_{T,s}(h) - C_s(\xi))|^2 \to 0,
\]
as $ T \to \infty $, provided $ N_T / T = o(1) $.
\end{theorem}

We shall now extend the above results to study weak consistency of the sequential cross-validation bandwidth under fairly general and weak assumptions. Having in mind the fact that $ h \sim T $, let us simplify the setting by strengthening that assumption to 
\begin{equation}
\label{HSpecific}
  h = h( \xi ) = T / \xi, \qquad \xi \in [1, \Xi],
\end{equation}
for some fixed $ \Xi \in (1,\infty) $. This means, $h$ and $ \xi $ are now equivalent parameters for each $T$. In what follows, we optimize over a compact interval, which is not restrictive for applications. Now $ \wh{m}_{h,-i} $ can be written as
\[
  \wh{m}_{h,-i} = \frac{1}{(i-1)h} \sum_{j=1}^{i-1} K( \xi (i-j)/T ) Y_j.
\]
With some abuse of notation, let us also write
\[
  C_{T,s}(\xi) = C_{T,s}( T/\xi ),
\]
i.e. from now on the expression $ C_{T,s}( T/\xi ) $ is studied as a function of $ \xi \in \Xi $.

The optimal cross-validated bandwidth is now given by $ h_T^*(s) = T/\xi^*_T(s) $, where 
\[
  \xi_T^*(s) = \argmin_{\xi \in \Xi} C_{T,s}(\xi),
\]
if $ C_{T,s} $ has a unique minimum; otherwise one selects a minimizer from the set $  \argmin_{\xi \in \Xi} C_{T,s}(\xi) $.

The next theorem yields weak consistency of the sequential cross-validation objective, uniformly over compact sets for the parameter $ \xi $ as well as uniformly over $ s \in \calS_N $, where again $ \calS_N $ is the set of time points (\ref{ST:SN}) where cross-validation is performed. 

\begin{theorem}  
\label{ST:Th4} Assume (\ref{ST:AssumptionsM}), (\ref{ST:Moments}), (\ref{ST:AssumptionsK}) and 
\eqref{ST:SN} such that $ N_T/T = o(1) $. Then, provided the bandwidth satisfies (\ref{HSpecific}), we have
\begin{equation}
\label{ST:Res1}
  \sup_{s \in \calS_N} \sup_{\xi \in [1,\Xi]} | C_{T,s}(\xi) - E C_{T,s}(\xi) | = o_P(1),
\end{equation}
and
\begin{equation}
\label{ST:Res2}
  \sup_{s \in \calS_N} \sup_{\xi \in [1,\Xi]} | C_{T,s}(\xi) - C_\xi(s) | = o_P(1),
\end{equation}
as $ T \to \infty $.
\end{theorem}

\begin{remark}
Notice that Theorem~\ref{ST:Th4} implies weak consistency of many other functionals, for example weighted loss functionals
\[
  \int_\Xi \int_S L( C_{T,s}(\xi) - C_\xi(s) ) w(s,\xi) ds \, d\xi
\]
where $ L $ is a Lipschitz continuous function attaining nonnegative values, $ S \subset [s_0,1] $ is a measurable set and $ w(s,\xi) $ is an integrable weighing function.
\end{remark}

\subsection{Consistency of the cross-validated bandwidth}

We are now in a position to formulate the following result on the asymptotic behavior of the cross-validated sequential bandwidth selector. The results of the previous subsection assuming i.i.d. errors are strong enough to apply known techniques to establish the consistency of minimizers of a sequence of random functions. 

\begin{theorem} 
\label{ST:Th5}
Suppose (\ref{ST:AssumptionsM}), (\ref{ST:Moments}), (\ref{ST:AssumptionsK}), 
\eqref{ST:SN} such that $ N_T/T = o(1) $ and (\ref{HSpecific}) hold true. Further, assume that one of the following 
conditions is satisfied,
\begin{itemize}
  \item[(i)] $ C_\xi(s) $ possesses a well-separated minimum $ \xi^* = \xi^*_s \in [1,\Xi] $, i.e.,
	\[
	  \inf_{\xi \in [1,\Xi] : |\xi - \xi^*| \ge \varepsilon} C_\xi(s) > C_{\xi^*}(s),
	\]
	for every $ \varepsilon > 0 $, or
  \item[(ii)] $ C_{T,s}( \xi ) $ is differentiable w.r.t. $ \xi $ such that $ \xi \mapsto \frac{\partial C_{T,s}(\xi)}{\partial \xi} $ is
  continuous and has exactly one zero.
\end{itemize}
Then
\[
  \xi_T^*(s) = \operatorname{argmin}_{\xi \in [1,\Xi]} C_{T,s}(\xi) \stackrel{P}{\to} \xi^*_s,
\]
as $ T \to \infty $.
\end{theorem}

Theorem~\ref{ST:Th5} asserts that the cross-validated bandwidth $ h^*_i $ computed at time $ t_i = s_i T $, $ s_i \in (s_0,1) $, is approximately given by $ \xi_{s_i}^* T $ for large $ T $. Notice that for given $ K $ and $ m $ the constant of proportionality, $ \xi_{s_i}^* $, can be calculated using the explicit formulas of Theorem~\ref{ST:Th0}.

\begin{remark}
There exist various sufficient criteria for consistency of argmin/argmax estimators, cf.  \citet{VaartAsymptoticStatistics}. Condition (i) is perhaps quite suited to the present problem, since our assumptions already ensure uniform convergence of the sequential cross-validation criterion and the requirement of a well-separated minimum can be checked analytically or numerically in an application for given $(K, \xi)$ and hypothesized $m$. Condition (i) ensures that the minimum is unique and especially rules out plateaus. Condition (ii) is a sufficient criterion, which is sometimes easier to verify, but requires the function $ \xi \mapsto C_\xi(s) $ to be differentiable with a continuous derivative.
\end{remark}

\section{Extensions to weakly dependent processes}
\label{Sec:Dependent}

Many series to which detection procedures are applied are dependent time series. This applies to almost all data sets arising in econometrics, environmetrics and communication engineering, but also to many data collected in biometrics and social sciences, e.g. longitudinal data in clinical trials or social surveys. In this case, procedures assuming i.i.d. error terms are not guaranteed to be valid. It is therefore quite natural to ask whether the results of the previous section carry over to the dependent case. Since often the specification of a parametric time series models for the error terms $ \{ \epsilon_t \} $ is subject of scientific discussion, we prefer to work with qualitative assumptions. 

Recall the definition of the $ \alpha $-mixing coefficient introduced by  \citet{Rosenblatt1956}. Let $ \{ Z_t \} $ be a weakly stationary process in discrete time. Let $ \calF^t = \sigma( Z_i : i \le t ) $ and $ \calF_{t} = \sigma( Z_i : i \ge t ) $. Then
\[
  \alpha(k) = \sup_{ A \in \calF^t, B \in \calF_{t+k} } | P(A \cap B) - P(A) P(B) |.
\]
The series $ \{ Z_t \} $ is called $ \alpha $-mixing if $ \lim_{k \to \infty} \alpha(k) = 0 $. For more information on mixing conditions we refer to \citet{Bosq1998}. Many commonly used (parametric) time series models are $ \alpha $-mixing. For instance, \citet{CarrascoChen2002} establish this property for ARCH models under certain conditions.

If $ \{ Z_t \} $ is $ \alpha $-mixing, series of the form $ g( Z_{t-m}, \dots, Z_{t+l}) $, $ m, l \in \N_0$, $g$ a measurable function, inherit that property. However, infinite functions of $ \alpha $-mixing processes are not necessarily mixing. As shown by \citet{Andrews84}, an $ AR(1) $ process with i.i.d. Bernoulli errors provides a well known counter-example. A more general notion is the following, which covers such processes.

\begin{definition}
\label{DefNED}
$ \{ Z_t \} $ is $L_r$-NED on $ \{ \xi_t \} $, $r> 0$, if there exist nonnegative constants $ \{ d_t : t \ge 1 \} $ with $ d_t \le 2 \| X_t - E(X_t) \|_r $ and $ \{ \nu_l : l \ge 0 \} $ such that 
\[
  \| Z_t - E( Z_t | \calF_{t-l}^{t+l} ) \|_r \le d_t \nu_l,
\]
and $ \nu_l \downarrow 0 $, as $ l \to \infty $, where $ \calF_s^t = \sigma( \xi_i : s \le i \le t ) $. 
\end{definition}

Compared to $ \alpha $-mixing, near epoch dependence can be viewed as a bridge to parametrically motivated models such as ARMA models or, more generally, linear processes. In our further discussion, we shall focus on $ L_2 $-NED. A $ L_2 $-NED series has the property that, by definition, one can approximate $ Z_t $ by its optimal $ L_2 $-predictor $ H( \xi_{t-l}, \dots, \xi_{t+l} ) = E( Z_t | \xi_{t-l}, \dots, \xi_{t+l} ) $ w.r.t. the $ L_2 $-norm, i.e. for any $ \varepsilon > 0 $ one can select $ l $ such that the $ L_2 $-error 
\[ 
  \nu_2 = \| Z_t - H( \xi_{t-l}, \dots, \xi_{t+l} ) \|_2
\]
does not exceed $ \varepsilon $. Parametrically motivated models are usually based on some i.i.d. noise process $ \{ \xi_t \} $. Let us suppose that for some function $g $ defined on $ \R^\infty $
\[
  Z_t = g( \dots, \xi_{t-1}, \xi_t, \xi_{t+1}, \dots ), \qquad t \in \Z,
\]
then
\[
  H( \xi_{t-l}, \dots, \xi_{t+l} ) = \int \cdots \int g( \dots, z_{-l-1}, \xi_{t-l}, \dots, \xi_{t+l}, z_{l+1}, \dots ) \, d \prod_{z_j : | j | > l} F(z_j),
\]
where $ F $ denotes the common d.f of the $ \xi_t $s. If, in addition, $g(z) = \sum_i \theta_i z_i $, $  z_i \in \R $ for $i \in \Z $, is a linear function with coefficients
$ \theta_i \in \R $ satisfying $ \sum_i | \theta_i | < \infty $, i.e. 
\[ 
  Z_t = \sum_i \theta_i \xi_{t-i}, \qquad t \in \Z
\] 
is a linear process, then we obtain
\[
  H( \xi_{t-l}, \dots, \xi_{t+l} ) = \sum_{i=t-l}^{t+l} \theta_i \xi_{t-i}, \qquad t \in \Z.
\]
The $L_2$-error, given by $ \nu_2 = 2( \sum_{|i| > l} | \theta_i | ) \| \xi_1 \|_2 $, converges to $ 0 $, if $ l \to \infty $, since the coefficients from a $l_1$ sequence.Thus, such linear processes are $ L_2 $-NED.

The following theorem provides the law of large numbers and weak consistency under strong mixing as well as a under a near epoch dependence condition.

\begin{theorem}
\label{ST:Th6} Let $K$ be a bounded kernel and suppose that
\begin{itemize}
  \item[(i)] $ \{ \epsilon_t \} $ is a weakly stationary $ \alpha $-mixing series with mixing coefficients $ \alpha(k) $, such that
$  \lim_{k \to \infty} k \alpha(k) = 0 $, $ T^{-1} \sum_{j,j'} | \Cov( \epsilon_j, \epsilon_{j'} ) | < \infty $ and $ T^{-2} \sum_{j,k,j',k'} | \Cov( \epsilon_j \epsilon_k, \epsilon_{j'} \epsilon_{k'} ) | < \infty $ hold true, or 
\item[(ii)]
$ \{ \epsilon_t \} $ is weakly stationary with $ E |\epsilon_1|^r < \infty $ for some $ r > 2 $ and $ L_2 $-NED on a weakly stationary $ \alpha $-mixing process such that the above conditions hold true for that underlying $\alpha$-mixing process. 
\end{itemize}
For fixed $ N \in \N $ put $ \calS_N = \{ s_1, \dots, s_N \} $ for given points $ 0 < s_0 \le s_1 < \cdots < s_N $. Then the uniform weak law of large numbers
holds true,
\[
  \max_{s \in \calS_N} \sup_{\xi \in [1,\Xi]} | C_{T,s} - C_\xi(s) |  = o_P(1),
\]
as $ T \to \infty $. Further, if $ C_\xi(s) $ satisfies condition (i) or (ii) of Theorem~\ref{ST:Th5}, then for fixed $ s \in [s_0,1] $
\[
  \operatorname{argmin}_{\xi \in [1,\Xi]} C_{T,s}(\xi) \stackrel{P}{\to} \xi^*_s,
\]
as $ T \to \infty $.
\end{theorem}

Notice that the results work under a less restrictive moment assumption and also allow for a larger class of kernels. Indeed, the kernel may have an unbounded support and is allowed to take negative values. However, we have to assume that the number $N$ of time points at which cross-validation is conducted is fixed. 


\section{Illustration: An Application in Photovoltaics}
\label{Sec:Illustration}

To illustrate the approach, we report about the following simulation experiment where we applied the method to a photovoltaic problem using real data to simulate error terms. In photovoltaics the power output of photovoltaic modules is the most important quantity for quality assessment, cf. \citet{StelandHerrmann2010}. In a scenario analysis we simulated a series of measurements according to the change-point model
\[
  Y_t  = \mu(t; \theta  ) + \epsilon_t, \qquad t = 1, \dots, T,
\]
where
\[
  \mu(t; \theta ) = \left\{
  \begin{array}{ll}
  \mu_0, & 1 \le t < q_1, \\
  \mu_0 + (t-q_1) \delta_1, & q_1 \le t < q_2, \\
  \mu_0 + (q_2-q_1) \delta_1 +  \Delta, & q_2 \le t,
  \end{array} \right.
\]
%
for $ t = 1, \dots, T $ with $ \theta = (\delta_1, \delta_2, q_1, q_2)' $ and $ T = 386 $.  $ \mu_0 = 200 $ denotes the nominal (target) power output, the parameters $ \delta_1 = -0.1 $ and $ \Delta $ model a drifting decreasing quality in terms of the mean power output with breaks (change-points) at $ q_1 = \trunc{T/4} $ and $ q_2 = \trunc{T/2} $. If $ \Delta = 0 $, then the process stabilizes after $ q_2 $ having a constant mean of $ 180.7 $; for this case study modules with power output larger than $ 180 $ were regarded as acceptable after re-labelling. Otherwise, there is a level shift of size $ \Delta $; for the scenario analysis we put $ \Delta = 2 s \approx 4.3 $ where $ s $ denotes the standard deviation of the errors. The errors were drawn from a real data set of photovoltaic measurements  to mimic a real-world setting in the simulation. $ \Delta = 0 $ represents the null model (null hypothesis $ H_0 $) of interest for our analysis.

\begin{center}
\begin{figure}
\includegraphics[width=9cm]{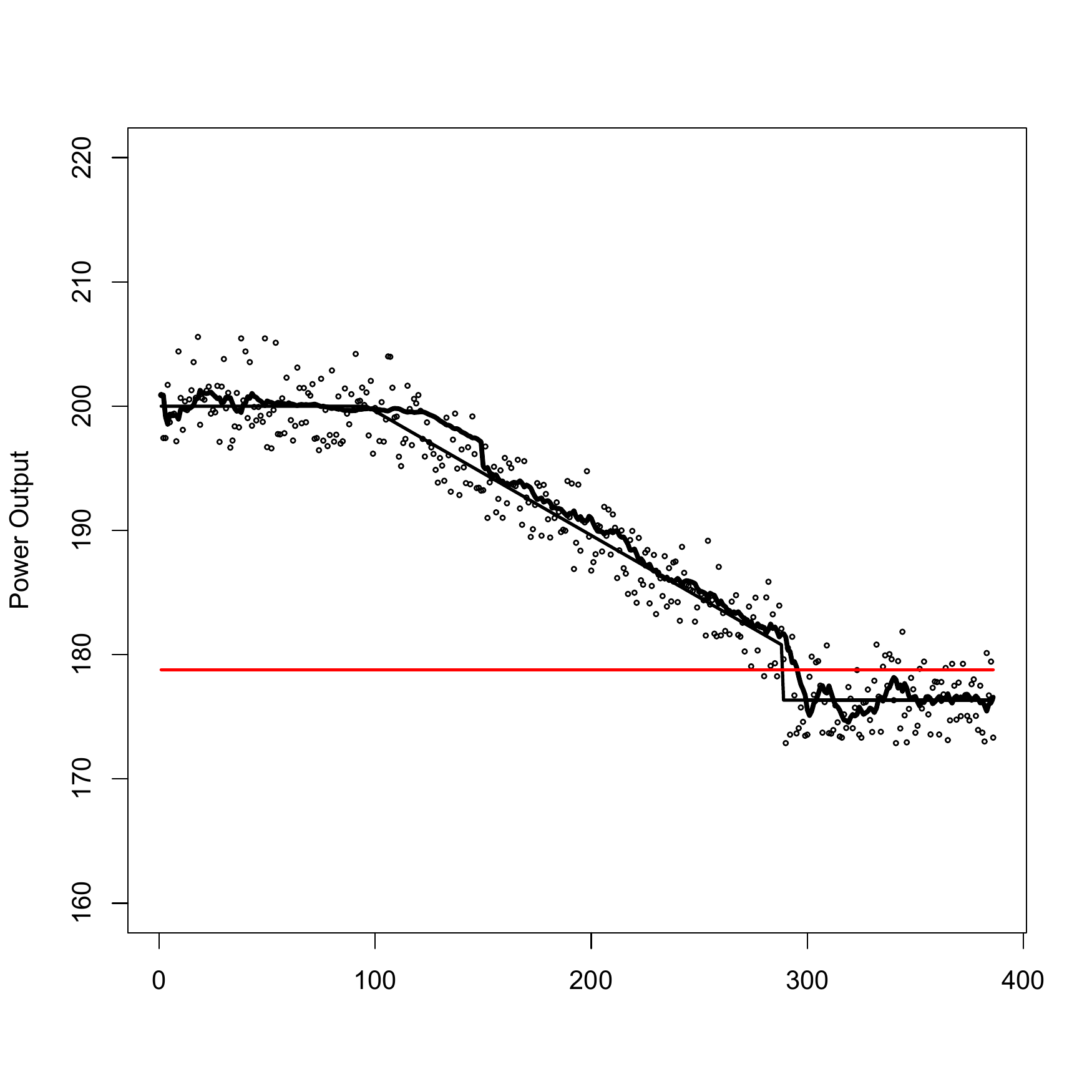}
\caption{A scenario analysis for photovoltaic measurements. True model (thin line), sequential smooth (bold) and control limit (red horizontal line).}
\label{Fig1}
\end{figure}
\end{center}

The simulated data and the cross-validated sequential kernel smooth are depicted in Figure~\ref{Fig1} and accompanied by a control limit $ c $. It can be seen that the predictions $ \widehat{m}_{h,-i} $ provide a reasonable approximation to the process mean. We applied the detector $ S_T^- $ using a Gaussian kernel, i.e., $ K(z) = (2\pi)^{-1/2}\exp( -z^2/2 ) $, 
$ z \in \R $, and the cross-validated bandwidth where cross-validation was conducted at the time points $ 50, 100, \dots, 350 $. The start of monitoring was determined using the rule $ \min(25, h_T^*(s_1) ) $ and the control limit was obtained by Monte
Carlo simulation ensuring an in-control average run length of $ ARL_0(S_T^-) = E_0(S_T^-) \approx 350 $ yielding $ c \approx 178.79 $. The signal is given at time instant $ 296 $. Table~\ref{Table1} provides simulation estimates of the mean delays defined for our purposes by $ E \max(0,S_T^- - q_2) $, as a function of $ \Delta $. Again,
the errors were simulated from real measurements. It can be seen that the chart reacts quickly to jumps.

\begin{table}
\caption{Simulation estimates of the mean delay for some values of $ \Delta $.}
\label{Table1}
\begin{center}
\begin{tabular}{l|ccccc} \hline
$\Delta$ &  $ 2/3 $ & $2/3$ & $ 4/3$ & $ 2 $ & $4$ \\ \hline
Mean delay & $60.25$ & $37.39$ & $10.01$ & $4.12$ & $2.07$ \\ \hline
\end{tabular}
\end{center}
\end{table}

\section*{Acknowledgments}

The author thanks Dr. W. Herrmann, T\"UV Rheinland Group (Cologne, Germany), for providing the real data and gratefully acknowledges financial support from DFG (German Research Foundation) and DAAD (German Academic Exchange Service), respectively. Remarks and hints from anonymous referees, which improved the presentation of the results, are also appreciated.

\appendix

\section{Proofs}


\begin{proof} (of Theorem~\ref{ST:Th0}).\\ 
Since the norming function $ s \mapsto N_T(s) = h^{-1} \sum_{i=1}^{\trunc{Ts}} K([\trunc{Ts}-i]/h) $, $ s \in 
[0,1] $,
is deterministic and converges to $ N_\xi(s) = \xi \int_0^s K(\xi(s-r)) \, d r $, we may and will assume that $ N_T(s) = 
N_\xi(s) = 1$.
First note that for $ j < i $ we have $ E( Y_i Y_j ) = m( i/T ) m( j/T ) $, since
$ E( \epsilon_j ) = E( \epsilon_i ) = 0 $ and $ E( \epsilon_i \epsilon_j ) = 0 $ by independence.
Further, $ E(Y_j^2) = m(j/T)^2 + E( \epsilon_j^2 ) $.
We have the decomposition $ E( C_{T,s}(h) ) = \wt{J}^{(1)}_{T,s} + \wt{J}^{(2)}_{T,s} + \wt{J}^{(3)}_{T,s} $ where
\begin{align*}
  \wt{J}^{(1)}_{T,s} & = -\frac2T \sum_{i=2}^{\trunc{Ts}} \frac{1}{h} \sum_{j=1}^{i-1} K([i-j]/h) E(Y_i Y_j), \\
  \wt{J}^{(2)}_{T,s} & = \frac1T \sum_{i=2}^{\trunc{Ts}} \left( \frac{1}{h} \right)^2 \sum_{j,k=1, j \not= k}^{i-1} K([i-j]/h)K([i-k]/h) 
E(Y_j Y_k), \\
  \wt{J}^{(3)}_{T,s} & = \frac1T \sum_{i=2}^{\trunc{Ts}} \left( \frac{1}{h} \right)^2 \sum_{j=1}^{i-1} K([i-j]/h)^2 E( Y_j^2 ).
\end{align*}
We provide the arguments for the second more involved term $ \wt{J}^{(2)}_{T,s}(h) $, the other terms are treated similarly. Notice that
\begin{align*}
  \wt{J}^{(2)}_{T,s} & =
   \frac1T \sum_{i=2}^{ \trunc{Ts} } \left( \frac{1}{h} \right)^2 
   \sum_{j,k=1, j \not= k}^{i-1} K([i-j]/h)K([i-k]/h) m(j/T) m(k/T) \\
  & = \int\limits_{\frac2T}^{ \frac{ \trunc{Ts} }{T} }  \left[ \frac{T}{h} \right]^2
   \int\limits_{\frac1T}^{ r - \frac1T }
   \int\limits_{\frac1T}^{ r - \frac1T }
   f_{T/h}(u,v,r) \, du dv dr
\end{align*}
where
\[
  f_{T/h}(u,v,r) =  K\left( \frac{ \trunc{Tr} - \trunc{Tu} }{h} \right) K\left( \frac{ \trunc{Tr} - \trunc{Tv} }{h} \right)
   m \left( \frac{ \trunc{Tu} }{T} \right) m\left( \frac{ \trunc{Tv} }{T} \right)
\]
Since $K$ as well as $ m $ are Lipschitz continuous and bounded, we have
\[
  f_{T/h}(u,v,r) \to f_\xi(u,v,r) = \xi^2 K( \xi(r-u) ) K( \xi (r-v) ) m( u ) m( v ),
\]
as $ T \to \infty $, uniformly in $ u,v,r \in [s_0,1] $. Notice that point-wise convergence holds true under weaker conditions, e.g. $ K$ bounded and continuous and $m$ continuous with $ \int_0^1 m^2(t) \, dt < \infty $, by dominated convergence. Now the result follows easily for the case $ N_T = 1 $, 
\[
  \wt{J}^{(2)}_{T,s} \to \wt{J}^{(2)}_s = \xi^2 \int_0^s \int_0^r \int_0^r K(\xi(r-u)) K(\xi(r-v)) m(u) m(v) \, du \, dv \, dr,
\]
uniformly, as $ T \to \infty $.
To  handle the general case, consider the analogous decomposition, $ C_{T,s}(h) = J^{(1)}_{T,s} + J^{(2)}_{T,s} + J^{(3)}_{T,s} $, where
for instance $ J^{(2)}_{T,s} = \wt{J}^{(2)}_{T,s} / N_T^2(s) $. Put $ J_s^{(2)} = \wt{J}^{(2)}_s / N_\xi^2(s) $ and notice 
that
\begin{equation}
\label{ST:Erw}
  \frac{ E( \wt{J}^{(2)}_{T,s}) }{ N_T^2(s) } - \frac{ \wt{J}^{(2)}_s }{ N_\xi^2(s) }
  =
  \frac{ N_\xi^2(s) E( \wt{J}^{(2)}_{T,s} ) - N_T^2(s) \wt{J}^{(2)}_s }{ N_\xi^2(s) N_T^2(s) }.
\end{equation}
Since $ K > 0 $ on $ (0,1) $,
\[
  \sup_{s \in [s_0,1]} N_T^{-2}(s), \ \sup_{s \in [s_0,1]} \left(
    \int_0^s K(\xi(s-u)) \, du \right)^{-2} = O(1),
\]
provided $ T $ is large enough. Further,
\[
  \sup_s | N_T^2(s) - N_\xi^2(s) |
  \le
  2 \sup_s | N_T(s) - N_\xi(s) | O( \| K \|_\infty \trunc{Ts}/h ) = o(1).
\]
Note that the numerator in \eqref{ST:Erw} equals
\[
  [ N_\xi^2(s) - N_T^2(s) ] \wt{J}^{(2)}_{T,s} + [ E( \wt{J}^{(2)}_{T,s} ) - \wt{J}^{(2)}_s ] N_\xi^2(s),
\]
which converges to $0$, uniformly in $ s \in [s_0,1] $, since $ \| \wt{J}^{(2)} \|_\infty,\ \| N_\xi \|_\infty < \infty $.
\end{proof}

\begin{proof} (of Theorem~\ref{ST:Th1})\\ 
To simplify the proof let us first assume that $ N_{T,-i} = 1 $.
Recall that for $j < i$ we have $ E( Y_i Y_j ) = m(i/T) m(j/T) $, since
$ E( \epsilon_j ) = E( \epsilon_i ) = 0 $ and $ E( \epsilon_i \epsilon_j ) = 0 $ by independence. Thus,
\begin{equation}
\label{DefZeta1}
  \zeta_{ij} = Y_i Y_j - E( Y_i Y_j ) = \epsilon_i m(j/T) + m(i/T) \epsilon_j + \epsilon_i \epsilon_j
\end{equation}
and
\begin{equation}
\label{DefZeta2}
  \zeta_{jj} = Y_j^2 - E( Y_j^2 ) = \epsilon_j^2 - E(\epsilon_j^2) + 2 \epsilon_j m(j/T).
\end{equation}
Notice that $ E( \zeta_{ij} ) = 0 $, and $ E( \zeta_{ij}^4 ) < \infty $, since, e.g.,
$ E( \epsilon_i \epsilon_j )^4 = ( E \epsilon_1^4 )^2 $ and
$ E( \epsilon_i m(j/T) (\epsilon_i \epsilon_j)^3 ) \le \| m \|_\infty E(\epsilon_1^4) E|\epsilon_1|^3 $.
Also note that $ E(Y_j^2) = m(j/T)^2 + E( \epsilon_j^2 ) $.
Consider the decomposition
\begin{equation}
\label{DecompCT}
  C_{T,s}(h) - E( C_{T,s}(h) ) = \wt{U}_{T,s} + \wt{V}_{T,s} + \wt{W}_{T,s}
\end{equation}
where
\begin{align}
\label{DefU}
 \wt{U}_{T,s} &=  -\frac2T \sum_{i=2}^{\trunc{Ts}} \frac{1}{h} \sum_{j=1}^{i-1} K([i-j]/h) \zeta_{ij}, \\
\label{DefV}
 \wt{V}_{T,s} &=
    \frac1T \sum_{i=2}^{\trunc{Ts}} \frac{1}{h^2} \sum_{j,k=1, j \not=k }^{i-1}
    K([i-j]/h) K([i-k]/h) \zeta_{jk}, \\
\label{DefW}
 \wt{W}_{T,s} &=
    \frac1T \sum_{i=2}^{\trunc{Ts}} \frac{1}{h^2} \sum_{j=1}^{i-1}
    K([i-j]/h)^2 \zeta_{jj}.
\end{align}
By virtue of Loeve's $ C_r $-inequality
\begin{equation}
\label{ST:Loeve}
  E|C_{T,s}(h) - E(C_{T,s}(h))|^2 \le 4( E | \wt{U}_{T,s} |^2 + E | \wt{V}_{T,s} |^2 ) + 2 E | \wt{W}_{T,s} |^2.
\end{equation}
Let us first consider  $ E | \wt{V}_{T,s} |^2 $ which can be written as
\begin{align}
\label{ST:Mom}
    & \frac{1}{T^2 h^4} \sum_{i,i'=2}^{\trunc{Ts}} \sum_{j,k=1, j \not=k}^{i-1} \sum_{j',k'=1, j' \not=k'}^{i'-1} K([i-j]/h) 
K([i-k]/h)
    \times \\ \nonumber
    & \qquad \quad K([i'-j']/h) K([i'-k']/h) E(\zeta_{jk} \zeta_{j'k'}).
\end{align}
Notice that $ E( \zeta_{jk} \zeta_{j'k'} ) $ equals
\[
  E( \epsilon_j m(k/T) + m(j/T) \epsilon_k + \epsilon_j \epsilon_k )
   ( \epsilon_{j'} m(k'/T) + m(j'/T) \epsilon_{k'} + \epsilon_{j'} \epsilon_{k'} )
\]
and vanishes, if $ \{ j,k \} \cap \{ j', k' \} = \emptyset $ by independence. Since $K$ has support $[-1,1] $,
the sums $ \sum_{j\not=k}$ and $ \sum_{j'\not=k'} $ concern only terms with $ |i-j|,|i-k| \le h $ and $ |i'-j'|, |i'-k'| \le h 
$.
Now consider the remaining non-vanishing cases. For the $ O(h^2) $ terms with $ j = j' \not= k' = k $, we have
\[
  E( \zeta_{jk} \zeta_{j'k'} ) = E(\epsilon_1 m(j/T) + m(k/T) \epsilon_2 + \epsilon_1 \epsilon_2)^2,
\]
which is non-negative and bounded in $j,k,T$. For the $ O(h^3)$ terms where $ j = j' $ and $ k \not= k' $ notice that
$ j = j' \not= k' $. Therefore,
\[
  E( \epsilon_j m(k/T) m(j'/T) \epsilon_{k'} ) \stackrel{j \not= k'}{=} 0, \qquad
  E( \epsilon_j m(k/T) \epsilon_{j'} \epsilon_{k'} ) \stackrel{j \not= k'} = 0,
\]
\[
  E( m(j/T) \epsilon_k \epsilon_{j'} m(k'/T) ) \stackrel{k \not= k'}{=} 0, \qquad
  E( m(j/T) \epsilon_k m(j'/T) \epsilon_{k'} ) \stackrel{k \not= k'}{=} 0,
\]
\[
  E( m(j/T) \epsilon_k \epsilon_{j'} \epsilon_{k'}) \stackrel{k' \not= k, k' \not= j'}{=} 0,
  E( \epsilon_j \epsilon_k m(j'/T) \epsilon_{k'} ) \stackrel{k' \not= k, k' \not= j'=j}{=} 0,
\]
and $ E( \epsilon_j \epsilon_k \epsilon_{j'} \epsilon_{k'} ) \stackrel{k' \not\in \{k,j'=j\}}{=} 0 $.
Thus, $ E( \zeta_{jk} \zeta_{j'k'} ) = E( \epsilon_1^2 m(k/T) m(k'/T) ) $ which is non-negative and finite. This shows that the 
contribution of
these terms is not larger than $ O( T^{-2}h^{-4} h^3 ) = O( T^{-3} )$. We may summarize that
there exists a constant $ c $ not depending on $ s $ such that \eqref{ST:Mom} is not larger than $ c T^{-3} $.
Consider now
\[
  | \wt{U}_{T,s} |^2 = \frac{4}{T^2h^2} \sum_{i,i'=2}^{\trunc{Ts}}
    \sum_{j=1}^{i-1} \sum_{j'=1}^{i'-1}
    K( [i-j]/h ) K( [i'-j']/h ) \zeta_{ij} \zeta_{i'j'}.
\]
If $ \{ i,j \} \cap \{ i',j'\} = \emptyset $, then $ E( \zeta_{ij} \zeta_{i'j'} ) = 0 $ by independence.
If $ i= i' $ and $ j \not= j' $ ($O(Th^2)$ terms), or $ i \not= i' $ and $ j = j' $ ($O(T^2h)$ terms),
or $ i = i' $ and $ j= j' $ ($O(Th)$ terms), again we have  $ 0 \le E( \zeta_{ij} \zeta_{i'j'} ) < c < \infty $ for some
constant $ c $ yielding
\[
  E | \wt{U}_{T,s} |^2 \le c \| K \|_\infty^2 T^{-1},
\]
where $ c $ does not depend on $ s \in [0,1] $.
We may conclude that the resulting upper bound for \eqref{ST:Loeve} is $ O(T^{-1}) $ for all $ s \in [0,1] $ yielding
\[
  \sup_{s \in [0,1]} E|C_{T,s}(h) - E(C_{T,s}(h))|^2 = O( T^{-1} ).
\]
The proof is now completed as follows.
\begin{align*}
  E \sup_{s \in \calS_N} | C_{T,s} - E( C_{T,s} ) |^2
    & \le E \sum_{s \in \calS_N} | C_{T,s} - E( C_{T,s} ) |^2 \\
    & \le | \calS_N | \sup_{s \in [0,1]} E | C_{T,s} - E( C_{T,s} ) |^2 \\
    & = O( N T^{-1} ) = O(T^{-1} ),
\end{align*}
as $ T \to \infty $. 

Let us now discuss the modifications when using $ \wh{m}_{n,-i} $ instead of $ \wt{m}_{n,-i} $.
Denote the corresponding decomposition of $ C_{T,s}(h) - E( C_{T,s}(h)) $ by $ U_{T,s} + V_{T,s} + W_{T,s} $, 
where the kernel weights $ K([j-i]/h) $ are replaced by $ K([j-i]/h) / \sum_{k=1}^i K([k-i]/h) $. 
We show how to treat $ U_{T,s} $. Notice that $ E | U_{T,s}(h) |^2 $ can be represented as
\begin{equation}
\label{IRepr}
  \int_{2/T}^{\trunc{Ts}} \int_{2/T}^{\trunc{Ts}}
  \int_{1/T}^{r-1/T} w_T(r) w_T(r') \int_{1/T}^{r-1/T} \int_{1/T}^{r'-1/T}
  G_{T,h}(u,v,r,r') \, dv \, du \, dr \, dr'
\end{equation}
where
\[
  0 \le G_{T,h}(u,v,r,r')
  =
  K( [\trunc{Tu}-\trunc{Tr}]/h )K( [\trunc{Tv}-\trunc{Tr'}]/h ) E( \zeta_{\trunc{Tr},\trunc{Tu}} \zeta_{\trunc{Tr'}, 
\trunc{Tv}} )
\]
and
\begin{align*}
  w_T(r)    & = 1 \biggr / \sum_{j=1}^{\trunc{Tr}-1} K( [\trunc{Tr}-j ]/h ) 
   \to w(r) = 1 \biggr/ \int_0^r K(\xi(r-z)) \, dz. 
\end{align*}
Here $ E( \zeta_{\trunc{Tr},\trunc{Tu}} \zeta_{\trunc{Tr'}, \trunc{Tv}} ) $ stands for the function
\[
  (r,u,r',v) \mapsto E(\zeta_{ij} \zeta_{i'j'}) \eins( r \in \calI_i, u \in \calI_j, r' \in \calI_{i'}, v \in \calI_{j'} ),
\]
where $ \calI_i = [i/T,(i+1)/T) $ for $ 0 \le i \le T $.
Since $ \gamma = \inf_{s \in [s_0,1]} \int_0^s K(\xi(s-z)) \, dz  > 0 $,
the elementary fact $ |1/x_n - 1/x| = | (x-x_n)/(x x_n)| $ yields
$ \sup_r | w_T(r) - w(r) | = O(1/T) $. Thus, if we replace in (\ref{IRepr}) the functions $ w_T(r) $ and $ w_T(r') $ by their
limits, by nonnegativity the difference can be bounded by 
\[
  \sup_r |w_T(r)w_T(r') - w(r)w(r') | 
    \int_{2/T}^{\trunc{Ts}} \int_{2/T}^{\trunc{Ts}}
  \int_{1/T}^{r-1/T} \int_{1/T}^{r-1/T} \int_{1/T}^{r'-1/T}
  G_{T,h}(u,v,r,r') \, dv \, du \, dr \, dr' 
\]
which is of the order $ O(1/T) $. To estimate the remaining term, namely
\[
    \int_{2/T}^{\trunc{Ts}} \int_{2/T}^{\trunc{Ts}}
  \int_{1/T}^{r-1/T} w(r) w(r') \int_{1/T}^{r-1/T} \int_{1/T}^{r'-1/T}
  G_{T,h}(u,v,r,r') \, dv \, du \, dr \, dr',
\]
notice that $ \sup_r w(r) \le 1/ \gamma < \infty $. Thus, the expression in the last display is not larger
than \[ \gamma^{-2} \int_{2/T}^{\trunc{Ts}} \int_{2/T}^{\trunc{Ts}}
  \int_{1/T}^{r-1/T} \int_{1/T}^{r-1/T} \int_{1/T}^{r'-1/T}
  G_{T,h}(u,v,r,r') \, dv \, du \, dr \, dr', \] 
which equals $ \gamma^{-2} | E \wt{U}_{T,s}(h) |^2 $.
\end{proof}

\begin{proof} (of Theorem~\ref{ST:Th2})\\ 
  By virtue of the method of proof used to establish Theorem~\ref{ST:Th1}, the result follows at once
  from  $ E \sup_{s \in \calS_N} | C_{T,s} - E( C_{T,s} ) |^2 = O( N_T T^{-1} ) = o(1) $.
\end{proof}

\begin{proof} (of Theorem~\ref{ST:Th3})\\ 
Define for $ s \in \calS_N $
$  \mathcal{X}_s = C_{T,s}(h) - E( C_{T,s}(h) ) $ and 
$  \mathcal{Y}_s = E(C_{T,s}(h) ) - C_s(\xi). $
Due to Theorem~\ref{ST:Th1}, we have
$
 E \sup_{s \in \calS_N} | \mathcal{X}_s | = o(1)
$
and Theorem~\ref{ST:Th0} yields
$
 E \sup_{s \in \calS_N} | \mathcal{Y}_s | = o(1).
$
Using the estimate
$
  | \mathcal{X}_s + \mathcal{Y}_s | \le \max( 2 | \mathcal{X}_s |, 2 | \mathcal{Y}_s | ),
$
which yields
$
  | \mathcal{X}_s + \mathcal{Y}_s |^2 \le 4 \max( | \mathcal{X}_s |^2, | \mathcal{Y}_s |^2 )
  \le 4 ( | \mathcal{X}_s |^2 + | \mathcal{Y}_s |^2 ),
$
we obtain
\[
  E \sup_{s \in \calS_N} | C_{T,s}(h) - C_s(\xi) |^2
  = E \sup_{s \in \calS_N} | \mathcal{X}_s + \mathcal{Y}_s | 
  \le 4 \sup_{s \in \calS_N} | \mathcal{X}_s |^2 + 4 \sup_{s \in \calS_N} | \mathcal{Y}_s |^2,
\]
which completes the proof.
\end{proof}

\begin{proof} (of Theorem~\ref{ST:Th4})\\
By Theorem~\ref{ST:Th0} it suffices to verify \eqref{ST:Res1}.
We make use of the decomposition for $ C_{T,s}(h) - E C_{T,s}(h) $ obtained above with the substitution
$ h = T/\xi $ and discuss the corresponding term $ V_{T,s}(\xi) $ in detail. The other terms are treated analogously and 
omitted. Fix some $ \xi' \in [1,\Xi] $ and $ \delta > 0 $. For brevity of notation, we will use the notation $ \sup_\xi = \sup_{\xi \in (\xi'-\delta,\xi'+\delta)} $ for the next steps and put
 \begin{equation}
 \label{DefWeights}
   w_{T,ijk} = \xi^2 K( \xi(i-j)/T ) K( \xi(i-k)/T ).
 \end{equation}
 The inequality $ \sup_x |f(x)| \le | \sup_x f(x) | + | \inf_x f(x) | $ yields
 \[
   \sup_\xi | V_{T,s}(\xi) | \le 2( V_{T,s}^{(1)}(\xi') + V_{T,s}^{(2)}(\xi') + V_{T,s}^{(3)}(\xi') )
 \]
 where
 \begin{align*}
   V_{T,s}^{(1)}(\xi') & =
     \left| \frac{1}{Th^2} \sum_{i,j,k}^{(s)} \sup_\xi w_{T,ijk}(\xi) \zeta_{jk}
         -  \frac{1}{Th^2} \sum_{i,j,k}^{(s)} E( \sup_\xi w_{T,ijk}(\xi) \zeta_{jk} ) \right|, \\
   V_{T,s}^{(2)}(\xi') & =
     \left| \frac{1}{Th^2} \sum_{i,j,k}^{(s)} \inf_\xi w_{T,ijk}(\xi) \zeta_{jk}
         -  \frac{1}{Th^2} \sum_{i,j,k}^{(s)} E( \inf_\xi w_{T,ijk}(\xi) \zeta_{jk} ) \right|, \\
   V_{T,s}^{(3)}(\xi') &=
        \frac{1}{Th^2} \sum_{i,j,k}^{(s)} E (\sup_\xi - \inf_\xi) w_{T,i,j,k}(\xi) \zeta_{jk}.
 \end{align*}
 Here $ \sum_{i,j,k}^{(s)} $ signifies $ \sum_{i=2}^{ \trunc{Ts} } \sum_{j,k=1,j \not=k}^{i-1} $.
 Since $ E|\sup_\xi w_{T,ijk}(\xi) \zeta_{jk} |^8 \le \| K \|_\infty^8 \Xi | \zeta_{jk} |^4 < \infty $,
 one can verify that for fixed $ \xi' $
 \[
   E\biggl( \sup_{s \in \calS_N} V_{T,s}^{(1)}(\xi') \biggr)^2 = o(1) \quad \text{and} \quad
   E\biggl( \sup_{s \in \calS_N} V_{T,s}^{(2)}(\xi') \biggr)^2 = o(1)
 \]
 using the same arguments as in the proofs of Theorem~\ref{ST:Th2} and \ref{ST:Th3}. 
 Let us now show that for fixed $ \xi' $ the term $ V_{T,s}^{(3)}(\xi') $ can be made arbitrary small. Indeed,
 by Lipschitz continuity of $K$, we may choose $ \delta > 0 $ small enough to ensure that
 \[
   \max_{1 \le j,k \le T} ( \sup_\xi - \inf_\xi ) \xi K(\xi(i-j)/T) \xi K( \xi(i-k)/T )
 \]
 is arbitrary small. Thus, by boundedness of $K$ and the dominated convergence theorem, for any $ \varepsilon > 0 $
the mean of each of the $ \sum_{i=1}^{ \trunc{Ts} } i(i-1) = O( \trunc{Ts}^3 ) $ summands of $ Th^2 V_{T,s}^{(3)}( \xi' ) $ (recall the definition
 (\ref{DefV})) is smaller than $ \varepsilon $, uniformly in $ i,j,k $, if $ \delta  > 0 $ is sufficiently small. Then
 \begin{align*}
   E \biggl| \sup_{s \in \calS_N} V_{T,s}^{(3)}( \xi' ) \biggr| 
   & \le 
     \sup_{s \in \calS_N} \frac{1}{Th^2} \sum_{i,j,k}^{(s)} E \max_{1 \le j,k \le T} ( \sup_\xi - \inf_\xi ) \xi K(\xi(i-j)/T) \xi K( \xi(i-k)/T ) 
     | \zeta_{jk} | \\
   & =
      \frac{1}{Th^2} \sum_{i,j,k}^{(1)} E \max_{1 \le j,k \le T} ( \sup_\xi - \inf_\xi ) \xi K(\xi(i-j)/T) \xi K( \xi(i-k)/T ) 
     | \zeta_{jk} | \\
   & = O( T^3 /(Th^2) \varepsilon ) = O( \varepsilon ).
 \end{align*}
 By compactness, a finite number of open balls
 $ B(\xi,\delta) = (\xi-\delta, \xi + \delta) $ cover $ [1,\Xi] $, such that 
 $ [1,\Xi] \subset \cup_{i=1}^M (\xi_i' - \delta, \xi_i' + \delta) $ for $ \xi_1', \dots, \xi_M' \in [1,\Xi] $, where
 $ M = M(\epsilon) \in \mathbb{N} $ depends on $ \varepsilon $.
 Thus, using the union bound we may now conclude that
 \begin{align*}
   P\biggl( \sup_{\xi \in [1,\Xi]} \sup_{s \in \calS_N}  | V_{T,s}(\xi) | > \varepsilon \biggr)
   & \le
   P\biggl( \max_{1 \le l \le N} \sup_{\xi \in B(\xi_l',\delta)}  \sup_{s \in \calS_N}  | V_{T,s}(\xi) |  > \varepsilon \biggr)  \\
   & \le
   \sum_{l=1}^M P\biggl( \sup_{s \in \calS_N} \sup_{\xi \in B(\xi_l',\delta)} V_{T,s}(\xi) > \varepsilon / M \biggr) \\
   & \le
    \sum_{l=1}^M 
   P\biggl( 2 \sup_{s \in \calS_N} \sup_{\xi \in B(\xi_l',\delta)}  \biggl| \sum_{\nu=1,2,3} V_{T,s}^{(\nu)}(\xi_l') \biggr| > \varepsilon / M \biggr) \\
   & \le
    \sum_{l=1}^M 
   P\biggl( 2 \sup_{s \in \calS_N}  \biggl| \sum_{\nu=1,2,3} V_{T,s}^{(\nu)}(\xi_l') \biggr| > \varepsilon / M \biggr).
\end{align*}   
We arrive at
\begin{equation}
\label{LLNEst}
   P\biggl( \sup_{\xi \in [1,\Xi]} \sup_{s \in \calS_N}  | V_{T,s}(\xi) | > \varepsilon \biggr)
   \le 
    \sum_{l=1}^M  \sum_{\nu=1,2,3} P\biggl( \sup_{s \in \calS_N} V_{T,s}^{(\nu)}( \xi_l' ) > \varepsilon /(6M) \biggr),
\end{equation}
with
\[
P\biggl( \sup_{s \in \calS_N} V_{T,s}^{(\nu)}( \xi_l' ) > \varepsilon /(6M) \biggr) =   O \left(
    E \biggl( \sup_{s \in \calS_N} V_{T,s}^{(\nu)}(\xi') \biggr)^2  \right), \qquad \nu = 1,2,
\]
and
\[
P\biggl( \sup_{s \in \calS_N} V_{T,s}^{(3)}( \xi_l' ) > \varepsilon /(6M) \biggr) =   O \left(
     E \biggl| \sup_{s \in \calS_N} V_{T,s}^{(3)}( \xi' ) \biggr| \right).
\]
This completes the proof, since $M$ is finite.
\end{proof}

\begin{proof} (of Theorem~\ref{ST:Th5})\\
  Confer \citet[Ch. 5.2]{VaartAsymptoticStatistics}.
\end{proof}

We shall make use of the following coupling lemma, cf. the works of \citet[Theorem~3]{Bradley82}, \citet{Schwarz80} and \citet[Lemma~1.2]{Bosq1998},
which allows to approximate directly dependent random variables by independent ones, a technique introduced in the papers 
\citet{BerkesPhillip77} and \citet{BerkesPhillip79}. 

\begin{lemma} (Bradley/Schwarz lemma).\\
\label{BradleysLemma}
Let $ (X,Y) $ be a $ \R^d \times \R $-valued random vector such that $ Y \in L_p $ for some
$ p \in [1,+\infty] $. Let $ c $ be a real number with $ \| Y + c \|_p > 0 $ and
$ \xi \in (0, \| Y + c \|_p ) $. Then there exists a random variable $ Y^* $ such that
$ P_{Y^*} = P_Y $, $ Y^* $ and $ Y $ are independent and
\[
  P(|Y-Y^*| > \xi ) \le 11(\xi^{-1} \| Y + c \|_p )^{p/(2p+1)} \bigl( \alpha(X,Y) \bigr)^{2p/(2p+1)}.
\]
Here $ \alpha(X,Y) $ denotes the $ \alpha $-mixing coefficient between $ \sigma(X) $ and $ \sigma(Y) $.
\end{lemma}

Notice that in Lemma~\ref{BradleysLemma} one may assume $ \alpha(X,Y) > 0 $, for otherwise the assertion
is trivially satisfied with $ Y^* = Y $, cf. \citet[p. 76]{Bradley82}. 

\begin{proof} (of Theorem~\ref{ST:Th6})\\
For weak consistency it suffices to show that either
\[
   \sup_{ \xi \in \Xi } | C_{T,s}(\xi) - C_\xi(s) | = o_P(1), 
\]
(and $ \xi \mapsto C_\xi(s) $ has a well-separated minimum), or
\[
  C_{T,s}(\xi) \to C_\xi(s), \qquad \text{for each $ \xi \in [1,\Xi] $}
\]
(and $ \xi \mapsto C_{T,s}(\xi) $ has an unique zero for each $ T $). Arguing along the lines of the proof of
Theorem~\ref{ST:Th4} we see that it suffices to show
\begin{equation}
\label{VTo0}
  V_{T,s}^{(k)}( \xi' ) \to 0, \qquad k = 1, 2,
\end{equation}
in probability. Indeed, since $ \calS_N $ is a finite set,  the right-hand side of (\ref{LLNEst}) is not larger than
the finite sum
\begin{equation}
\label{FiniteSum}
  \sum_{l=1}^M \sum_{\nu = 1, 2, 3} \sum_{s=1}^N P( V_{T,s}^{(\nu)} > \varepsilon / (6MN) ).
\end{equation}
For i.i.d. errors with fourth moments this is a consequence of our results on uniform $L_2$ convergence. 
To extend the result to the dependent case, we apply the Bradley/Schwarz lemma \ref{BradleysLemma} and the block-splitting
technique to the random variables $ V_{T,s}^{(k)}( \xi' ) $, $ k = 1, 2 $.  We have to check that their second moments
are bounded by some constant $ U < \infty $ (which depends on $ \| m \|_\infty $, $ \| K \|_\infty $, $ E(\epsilon_1^4) $ and $ \Xi $). 
Let us sketch the arguments for $ k = 1 $. Recall the definition (\ref{DefWeights}) and the representation
(\ref{ST:Mom}) to see that 
\begin{align} \nonumber
  E ( V_{T,s}^{(1)}(\xi') )^2 & \le \frac{1}{T^2h^4} \sum_{i,i'=2}^{\trunc{Ts}} \sum_{j,k=1, j \not=k }^{i-1} \sum_{j', k' = 1, j' \not=k' }^{i'-1}
  |w_{T,ijk} w_{T,i'j'k'} | | \Cov( \zeta_{jk}, \zeta_{j'k'} ) | \\
  \label{SecondMomentDep}
  & = O\left( \frac{1}{T^4}  \sum_{j,k=1, j \not=k }^{i-1} \sum_{j', k' = 1, j' \not=k' }^{i'-1} | Cov( \zeta_{jk}, \zeta_{j'k'} ) | \right).
\end{align}
Notice that the estimate (\ref{SecondMomentDep}) requires the kernel $K$ only to be bounded. 
Recalling $ \| m \|_\infty < \infty  $ and the definition of the $ \zeta_{jk} $ in (\ref{DefZeta1}) and (\ref{DefZeta2}), respectively, we get
\[
  | Cov( \zeta_{jk}, \zeta_{j'k'} ) | \le \max\{ \| m \|_\infty, 1 \} 
  ( | \Cov( \epsilon_j, \epsilon_{j'} ) | + \cdots + | \Cov( \epsilon_j \epsilon_k, \epsilon_{j'} \epsilon_{k'} ) | ).
\]
Estimating separately the resulting sums obtained when combining that inequality with (\ref{SecondMomentDep}), we see that $ E ( V_{T,s}^{(1)}(\xi') )^2 = O(1) $,
uniformly in $ s $, since $ T^{-1} \sum_{j,j'} | \Cov( \epsilon_j, \epsilon_{j'} ) | < \infty $ and $ T^{-2} \sum_{j,k,j',k'} | \Cov( \epsilon_j \epsilon_k, \epsilon_{j'} \epsilon_{k'} ) | < \infty $. In other words, the second moment of $ Th^2  V_{T,s}^{(1)}(\xi')  $, a sum of $ O(T^3) $ summands, is $ O(T^6) $, and the same applies whenever selecting, say, $p$ summands and considering the second moment of their sum $S_p$, i.e. $ E(S_p^2) = O(p^2) $.

Let  $ V_{T,s}(\xi') \in \{ V_{T,s}^{(k)}( \xi' ) : k = 1,2 \} $ and denote the $ m = O(Th^2) $ summands by
$ \rho_{ijk} $ such that $ m V_{T,s}(\xi') = \sum_{i=1}^{\trunc{Ts}} \sum_{j,k=1, j \not= k}^{i-1} \rho_{ijk} $. Let us now apply the block-splitting technique in combination
with the Bradley/Schwarz coupling lemma. Notice that the following derivations do not depend on the kernel at all.
For simplicity, we shall assume $ m = 2pq $, where $p$ and $q$ will be chosen later. Partition the $m$ summands in consecutive blocks of length $p$. Clearly, $m V_{T,s}( \xi' ) $ is the sum of the $ 2q $ partial sums of these blocks. Number these partial sums from $1$ to $2q$ and denote the partial sums corresponding to odd numbers by $ B_1, \dots, B_q $ and those corresponding to even numbers by $ B_1', \dots, B_q' $. Let $ \varepsilon > 0 $. It suffices to establish a bound for $ P( |B_j| > m \varepsilon ) $. Put  $ c = 2pU $ and notice that
\[
  \min_j \| B_j + c \|_2 \ge c - \max_j \| B_j \|_2 = pU,
\]
since $ \| B_j \|_2 \le p U $ by Minkowski's inequality, for all $j$. Let 
\[
  \xi = \min( pU, m \varepsilon / (4q) ) \in (0, \min_j \| B_j + c \|_2 ).
\]
Applying Bradley's lemma 
yields the existence of $ B_j^* $ independent from $ B_1^*, \dots, B_{j-1}^* $ with $ B_j^* \stackrel{d}{=} B_j $ and
\begin{align*}
  P( | B_j - B_j^* | > \xi ) & \le 11 \left( \frac{\|B_j+c\|}{\xi} \right)^{2/5} \alpha(p)^{4/5} \\
    & \le 11 \left( \frac{3pU}{\min( pU, m \varepsilon/(4q) ) } \right)^{2/5} \alpha(p)^{4/5} \\
    & \le 11(3+6U/\varepsilon)^{2/5} \alpha(p)^{4/5}.
\end{align*}
By independence of $ B_1^*, \dots, B_q^* $ and since $ E(B_j^*)^2 = O(p^2) $, we have
\[
  P\biggl( \biggl| \sum_{j=1}^q B_j^* \biggr| > m \varepsilon/4 \biggr) 
  \le \frac{ \sum_{j=1}^q E(B_j^*)^2 }{ m^2 \varepsilon^2/4 } 
  = O( 4/(2q\varepsilon^2) ).
\]
The inclusion
\[
\biggl\{ \biggl| 
  \sum_{j=1}^q B_j \biggr| > \frac{ m \varepsilon }{ 2 } \biggr\}
  \subset
  \biggl\{ \biggl| \sum_{j=1}^q B_j \biggr| > \frac{ m \varepsilon }{ 2 }, | B_j - B_j^* | \le \xi, 1 \le j \le q \biggr\}
  \cup
  \bigcup_{j=1}^q \{ | B_j - B_j^* | > \xi \}.
\]
leads us to 
\begin{align*}
P \biggl( \biggl| \sum_{j=1}^q B_j \biggr| > \frac{m\varepsilon}{2} \biggr)
& \le
  P \biggl( \biggl| \sum_{j=1}^q B_j^* \biggr| > \frac{m\varepsilon}{4} \biggr) +
   \sum_{j=1}^q P(|B_j-B_j^*| > \xi ) \\
& \le O( (q\varepsilon)^{-1} ) + O( q(3+6U/\varepsilon)^{2/5} \alpha(p)^{4/5} ).
\end{align*}
Putting $ p = \trunc{ m^{1/2} } $ and $ q = m/(2p) $, we may conclude that the last expression is $ o(1) $, provided
$ \lim_{k \to \infty} k \alpha(k) = 0 $. This shows that (\ref{VTo0}) holds for the strong mixing case, 
\begin{equation}
\label{VTo0Mixing}
  P( V_{t,s}^{(k)}( \xi' ) > \varepsilon ) \le P\biggl( \biggl| \sum_{j=1}^q B_j \biggr| >  \frac{m \varepsilon}{2} \biggr) 
    + P \biggl( \biggl| \sum_{j=1}^q B_j' \biggr| >   \frac{m \varepsilon}{2} \biggr) = o(1).
\end{equation}
Let us now discuss the case that $ \{ \epsilon_t \} $ is $ L_2 $-NED on a stationary $ \alpha $-mixing
process
such the the above arguments hold true when the error process $ \{ \epsilon_t \} $ is replaced by that $ \alpha $-mixing process. As in Definition~\ref{DefNED},
denote the corresponding $ \sigma $-fields by $ \calF_a^b $. Put $ \wt{\epsilon}_t = E(\epsilon_t | \calF_{t-l}^{t+l} ) $ 
for all $ t $ and $ l > 1$. Then $ \| \epsilon_t - \wt{\epsilon}_t \|_2 = O( \nu_l ) $ for some sequence
$ \nu_l = o(1) $, as $ l \to \infty $, uniformly in $t$. Recall the definitions (\ref{DefZeta1}) and (\ref{DefZeta2}). 
Denote by $ \wt{V}_{T,s}(h) $ the statistic $ V_{T,s}(h) $ where the $ \epsilon_i $ are replaced by the $ \wt{\epsilon}_i $ i.e.
$ \zeta_{ij} $ are replaced by the random variables
\begin{align*}
  \wt{\zeta}_{ij}
    = \wt{\epsilon}_i \wt{\epsilon}_j - E(\wt{\epsilon}_i \wt{\epsilon}_j ) + \epsilon_i m(j/T) + \epsilon_j m(i/T).
\end{align*}
We have 
\begin{align*}
  \| \zeta_{ij} - \wt{\zeta}_{ij} \|_1
    & \le
     \| \epsilon_i \epsilon_j - \wt{\epsilon}_i \wt{\epsilon}_j \|_1 + \| E( \epsilon_i \epsilon_j - \wt{\epsilon}_i \wt{\epsilon}_j ) \|_1 \\
     & \qquad + \| m(i/T)( \epsilon_j - \wt{\epsilon}_j ) \|_1 + \| m(j/T)( \epsilon_i - \wt{\epsilon}_i ) \|_1\\
     &\le
     2 \| \epsilon_i \epsilon_j - \wt{\epsilon}_i \wt{\epsilon}_j \|_1 + 2 \| m \|_\infty \| \epsilon_i - \wt{\epsilon}_i \|_1 \\
     &\le
     2\bigl( \| \epsilon_i \|_2 \| \epsilon_j - \wt{\epsilon}_j \|_2 + \| \wt{\epsilon}_j \|_2 \| \epsilon_i - \wt{\epsilon}_i \|_2  + \|m\|_\infty \| \epsilon_i - \wt{\epsilon}_i \|_2 \bigr) 
     = O( \nu_l ).
\end{align*}
For fixed $ \xi' $ we obtain
\begin{align*}
  \sup_{s \in [s_0,1]} \| V_{T,s}(\xi') - \wt{V}_{T,s}(\xi') \|_1
   & = \sup_{s \in [s_0,1]} \left\| \sum_{i=1}^{\trunc{Ts}} \sum_{j,k=1, j \not= k}^{i-1} w_{T,ijk}( \zeta_{ij} - \wt{\zeta}_{ij} ) \right\|_1 
     = O(\nu_l ),
\end{align*}
which allows us to estimate the summands in (\ref{FiniteSum}) by
\begin{align*}
  P( V_{T,s}^{(\nu)} > \varepsilon / (6MN) )
   & = P( \wt{V}_{T,s}^{(\nu)} > \varepsilon / (6MN) ) + \sup_{s \in [s_0,1]} \|  V_{T,s}^{(\nu)} -  \wt{V}_{T,s}^{(\nu)} \|_1 \\
   & \le  o(1) + O( \nu_l ),
\end{align*}
since $ \wt{V}_{T,s}^{(\nu)} $ is calculated from the $ \alpha $-mixing random variables $ \wt{\epsilon}_i $, such that
(\ref{VTo0Mixing}) applies. This completes the proof.
\end{proof}

\bibliographystyle{agsm}
\bibliography{lit}

\end{document}